\documentclass[10pt,runningheads,a4paper]{article}
\usepackage{latexsym}
\usepackage{epsfig}
\usepackage{epstopdf}
\usepackage{amssymb}
\usepackage{fancyhdr}

\usepackage{graphicx}
\usepackage{amsmath}
\usepackage{amssymb}
\usepackage{cite}

\makeatletter

\newcommand{\Rmnum}[1]{\expandafter\@slowromancap\romannumeral #1@}
\makeatother

\usepackage{amsmath,amssymb,graphics,float}
\usepackage{pstricks,pst-node,pst-tree,pst-plot}
\usepackage{graphicx,epsfig,subfigure}

\title{  First Integrals of Dynamical Systems And Their Numerical Preservation }
\usepackage{authblk}
\author[1]{W. Irshad}
\author[2]{Y. Habib}
\author[3]{M. U. Farooq}
\affil[1]{School of Natural Sciences, National University of
Sciences and Technology, Sector H-12, Islamabad, Pakistan. }
\affil[2]{Department of Mathematics, COMSATS Institute of Information Technology, Lahore, Pakistan.  }
\affil[3]{Department of BS and H, College of E and ME, National
University of Sciences and Technology, Peshawar Road, Rawalpindi,
Pakistan.\newline Emails: wajeehairshad97@gmail.com,
yhabib@ciitlahore.edu.pk, m\textunderscore ufarooq@yahoo.com }
\date{}

\begin{document}
\maketitle
\begin{abstract}
\noindent We calculate Noether like operators and first integrals of scalar equation $y''=-k^2 y$ using complex Lie symmetry method, by taking values of $k$ and $y$ to be real as well as complex. We numerically integrate the equations using a symplectic Runge-Kutta method and check for preservation of these first integrals. It is seen that these structure preserving numerical methods provide qualitatively correct numerical results and good preservation of first integrals is obtained.
\end{abstract}
\vspace{.3cm}
\vspace{.3cm}
\section{Introduction}
Marius Sophus Lie proposed a symmetry
based method for the analytical solution of differential equations
using group of continuous transformations known as Lie groups
\cite{SL1, SL2, Ibragimov, Stephani}. Amalie Emmy Noether later presented her remarkable theorem which relates variational symmetries with conservation laws or the first integrals in \cite{EN}. In literature, different methods are available to calculate first
integrals of ordinary differential equations (ODEs) including the
direct method, the characteristic or multiplier method, the Noether
approach and partial Noether approach \cite{Naz, Leach, MZ, Kara }.
In this paper, we use classical Noether approach to calculate first integrals of harmonic oscillator equation. We then apply complex
symmetry method in the restricted domain to find first
integrals of system of harmonic oscillators by considering the
Lagrangian in the complex variable domain \cite{UF,AQ,SA}.
\\
\noindent Concerning the numerical solutions of ODEs with quadratic
first integrals, it is well known that symplectic numerical methods
are a suitable candidate \cite{Serna}. These methods are a subclass
of geometric integrators which preserve geometric properties of the
exact flow of ODEs. One class of symplectic methods with optimal order are the Gauss-Legendre Runge--Kutta methods. They are one step numerical methods for ODEs and preserve all linear and quadratic first integrals of a dynamical system \cite{rf4}. If we intend to preserve cubic or higher order first integrals, we do not have a general numerical scheme for such purpose but we can design a
numerical method that has this as its specific goal, for example splitting method and discrete gradient method \cite{rf4}. In this paper, we present a way of constructing symplectic Runge--Kutta methods. We then take order four Gauss-Legendre Runge--Kutta method for the numerical integration of ODEs
and report good preservation of first integrals by the numerical solution. 

\section{\bf Symmetries and First Integrals}\label{sfi}
Consider a second-order ODE,
\begin{equation}\label{1}
y''=f(t,y,y'),
\end{equation}
which admits a Lagrangian $L(t,y,y')$ satisfying the Euler-Lagrange
equation,
\begin{equation}\label{2}
\frac{\partial{L}}{\partial{y}}-\frac{d}{d{t}}(\frac{\partial{L}}{\partial{y{'}}})=
0.
\end{equation}
To explain the invariance criteria for variational problems under a group of transformation we consider the operator,
\begin{equation}\label{3}
X=\xi(t,y)\frac{\partial}{\partial{t}}+\eta(t,y)\frac{\partial}{\partial{y}},
\end{equation}
known as the infinitesimal operator or the symmetry generator. The
functions $\xi$ and $\eta$ are the components of tangent vector
${X}$ and are defined as,
\begin{equation}\label{11}
\xi(t,y)=\frac{\partial{\tilde{t}}}{\partial{\epsilon}}\mid_{\epsilon=0},
~~~~~~\eta(t,y)=\frac{\partial{\tilde{y}}}{\partial{\epsilon}}\mid_{\epsilon=0}.
\end{equation}
The operator $X$ is called Noether symmetry generator
corresponding to the Lagrangian $L(t,y,y')$, if there exist a gauge
function $B(t,y)$ such that the following condition holds,
\begin{equation}\label{4}
{X^{(1)}}(L)+D_{t}(\xi)L= D_{t}(B),
\end{equation}\\
where $X^{(1)}$ is a first order prolongation of $X$ and $D$ is
total derivative operator given by,
\begin{equation}\label{5}
D_{t}=\frac{\partial}{\partial{t}}+y'\frac{\partial}{\partial{y}}.
\end{equation}
According to Noether theorem, for each Noether symmetry of an
Euler-Lagrange equation, there corresponds a function $I$
\begin{equation}\label{6}
I=\xi{L}+(\eta-\xi y')\frac{\partial{L}}{\partial{y'}}-B(t,y),
\end{equation}
called the first integral or conserved quantity of the equation
\eqref{1}, with respect to the symmetry generator $X$.
\subsection{Complex Symmetry Analysis}
We first discuss some important results related to
complex Noether symmetries, complex Lagrangian and Noether Theorem
in the restricted complex domain. We will use them to determine 
first integrals of second-order restricted complex ODEs \cite{SAQ}.
We then present expressions for Euler-Lagrange like equations, conditions
for Noether-like operators and expressions for first
integrals
corresponding to these operators. For more details, see \cite{UF} and references therein.\\
Consider a system of two second-order ODEs of the form,
\begin{eqnarray}\label{7}
\begin{aligned}
f'' &= w_{1}(t,f,g,f',g'),\\
g'' &= w_{2}(t,f,g,f',g').
\end{aligned}
\end{eqnarray}\\
Suppose we have a transformation $y(t)=f+i g$ and $w=w_{1}+iw_{2}$
which converts the system \eqref{7} to second order restricted
complex ODE,
\begin{equation}\label{8}
y''=w(t,y,y').
\end{equation}\\
Assume that the equation \label{8} admits a complex Lagrangian
$L(t,f,g,f',g')$, i.e., $L=L_{1}+i L_{2}$. Therefore, we have two
Lagrangians $L_{1}$ and $L_{2}$ for the system \eqref{7} that
satisfy Euler-Lagrange like equations,
\begin{eqnarray}\label{9}
\begin{aligned}
\frac{\partial L_{1}}{\partial f} + \frac{\partial L_{2}}{\partial g}-\frac{d}{dt}(\frac{\partial L_{1}}{\partial f'} + \frac{\partial L_{1}}{\partial g'})&=0,\\
\frac{\partial L_{2}}{\partial f} - \frac{\partial L_{1}}{\partial
g}-\frac{d}{dt}(\frac{\partial L_{2}}{\partial f'} - \frac{\partial
L_{1}}{\partial g'})&=0.
\end{aligned}
\end{eqnarray}
The operators,
\begin{eqnarray}\label{10}
\begin{aligned}
X_{1}&=\varsigma_{1}\frac{\partial}{\partial{t}}+\chi_{1}\frac{\partial}{\partial{f}}+\chi_{2}\frac{\partial}{\partial{g}},\\ X_{2}&=\varsigma_{2}\frac{\partial}{\partial{t}}+\chi_{2}\frac{\partial}{\partial{f}}-\chi_{1}\frac{\partial}{\partial{g}}.
\end{aligned}
\end{eqnarray}
are called Noether-like operators for the Lagrangians $L_{1}$ and $L_{2}$, if they satisfy following conditions,
\begin{equation}\label{11}
\begin{aligned}
X^{(1)}_{1}(L_{1})-X^{(1)}_{2}(L_{2})+(D_{t}\varsigma_{1})L_{1}-(D_{t}\varsigma_{2})L_{2}&=D_{t}A_{1},\\
X^{(1)}_{1}(L_{2})+X^{(1)}_{2}(L_{1})+(D_{t}\varsigma_{1})L_{2}+(D_{t}\varsigma_{2})L_{1}&=D_{t}A_{2},
\end{aligned}
\end{equation}\\
where $A_{1}$ and $A_{2}$ are suitable gauge functions. The two
first integral corresponding to the Noether-like operators $X_{1}$
and $X_{2}$ can be found as,
\begin{equation}\label{12}
\begin{aligned}
I_{1}&=\varsigma_{1}L_{1}-\varsigma_{2}L_{2}+{\partial_{f'}}L_{1}(\chi_{1}-f'\varsigma_{1}-g'\varsigma_{2})-{\partial_{f'}}L_{2}(\chi_{2}-f'\varsigma_{2}-g'\varsigma_{1})-A_{1},\\
I_{2}&=\varsigma_{1}L_{2}+\varsigma_{2}L_{1}+{\partial_{f'}}L_{2}(\chi_{1}-f'\varsigma_{1}-g'\varsigma_{2})+{\partial_{f'}}L_{1}(\chi_{2}-f'\varsigma_{2}-g'\varsigma_{1})-A_{2}.
\end{aligned}
\end{equation}\\
\section{Runge--Kutta methods}\label{rkm}
\noindent Runge--Kutta methods \cite{rf3} are one-step numerical methods for solving initial value problems (IVPs),
\begin{equation}\label{ivp}
y'(t)=f(y(t)), \hspace{0.3in} y(t_0)=y_0, \hspace{0.3in} y(t)\in \mathbb{R}^{n}.
\end{equation}
These methods provide an approximation $y_{n}=y(t_{n})$ of the exact
solution $y(t)$ at time $t_{n} = nh$,
where $n = 0,1,\cdots$ and $h$ corresponds to the stepsize. The
general form of an $s$-stage Runge-Kutta method is,
\begin{align}\label{RK_general_form}
Y_{i} &= y_{n-1} + \displaystyle\sum_{j=1}^{s}a_{ij}hf(Y_{j}), \hspace{0.2in} i = 1,2,\cdots,s, \\ \nonumber
y_{n} &= y_{n-1} + \displaystyle\sum_{i=1}^{s}b_{i}hf(Y_{i}),
\end{align}
where $b_{i}$ are the quadrature weights of the method and $c_{i}$
are the nodes at which the stages $Y_{i}$ are evaluated. A
Runge--Kutta method can be represented by a Butcher tableau,
\begin{align*}
\begin{array}{c|cccc}
c_1 & a_{11} & a_{12} & \cdots & a_{1s}\\
c_2 & a_{21} & a_{22} & \cdots & a_{2s} \\
\vdots & \vdots & \vdots & \ddots & \vdots \\
c_s & a_{s1} & a_{s2} & \cdots & a_{ss} \\
\hline
& b_1 & b_2 & \cdots & b_s
\end{array}.
\end{align*}
The Runge--Kutta methods are explicit if $a_{ij} = 0$ for $i \leq j,
$ otherwise, they are implicit.
\subsection{Symplectic Runge--Kutta methods}
If the initial value problem \eqref{ivp}
has a quadratic first integral
\[I(y)= \langle y,S y \rangle = y^{T}Sy,\]
where $S$ is a symmetric square matrix, then we have
\[\langle y,f(y) \rangle=y^{T}Sf(y) = 0.\]
We want to determine numerical solutions $y_{n}$ such that the first
integral $I(y)$ is preserved numerically, i.e.,
\[
\langle y_{n},S y_{n} \rangle = \langle y_{n-1},S y_{n-1} \rangle
\hspace{0.2in} n=0,1,\hdots .\] It has been shown in
\cite{Stability_cooper,SRK_LAS,BJ} that only symplectic Runge--Kutta
methods preserve the quadratic first integrals while numerically integrating
\eqref{ivp}. Moreover, in this paper we will only be
considering implicit Runge--Kutta methods to check the numerical
preservation of first integrals because explicit methods cannot be
symplectic \cite{Sanz}. A Runge--Kutta method is
symplectic if its coefficients satisfy the following condition \cite{Stability_cooper,SRK_LAS,rf20},
\begin{equation}\label{V1}
b_{i}a_{ij}+b_{j}a_{ji}-b_{i}b_{j}=0 \hspace{.2cm}
\text{for\hspace{.2cm}all} \hspace{.3cm}i,j=1,2,\hdots,s,
\end{equation}
which can be derived as follows.\\ Firstly, apply the
Runge-Kutta method \eqref{RK_general_form} to solve the
IVP \eqref{ivp}. The stage values are
\[
Y_{i}=y_{n-1}+h \displaystyle \sum_{j=1}^{s}a_{ij}f(Y_{j}).
\]
Since,
\begin{align}\label{SRKS}
&\langle Y_{i},S f(Y_{i})\rangle = 0,\nonumber\\
\Rightarrow~ &\langle y_{n-1},S f(Y_{i})\rangle+h \displaystyle
\sum_{j=1}^{s}a_{ij}\langle f(Y_{j}),S f(Y_{i})\rangle=0.
\end{align}
Moreover, for the output values we have,
\[
y_{n} = y_{n-1} + \displaystyle\sum_{i=1}^{s}b_{i}hf(Y_{i}).
\]
Thus
\begin{align}\label{SRKI}
\langle y_{n},S y_{n}\rangle &=\langle y_{n-1},S y_{n-1}\rangle +h\displaystyle \sum_{i=1}^{s}b_{i}\langle y_{n-1},S f(Y_{i})\rangle \nonumber\\
&+h\displaystyle \sum_{j=1}^{s}b_{j}\langle f(Y_{j}),S y_{n-1}\rangle+h^{2}\displaystyle \sum_{i,j=1}^{s}b_{i}b_{j}\langle f(Y_{i}),S f(Y_{j})\rangle.
\end{align}
Evidently from \eqref{SRKS} and \eqref{SRKI} we have,
\[\langle y_{n},S y_{n}\rangle = \langle y_{n-1},S y_{n-1}\rangle ,\]
provided
\begin{equation}\label{sympl_id}
b_{i}a_{ij}+b_{j}a_{ji}-b_{i}b_{j}=0.
\end{equation}
\subsection{Construction of symplectic Runge-Kutta methods}
\label{sec:3}
Although there exist several techniques to construct symplectic Runge-Kutta methods in literature \cite{rf4,rf30}, here we construct symplectic Runge-Kutta methods with the help of Vandermonde transformation. This was first discussed in \cite{rf7}.  The idea is to pre and post multiply the Vandermonde matrix with the matrix of symplectic condition for Runge--Kutta method \eqref{V1}. Our strategy is to write the values of $a$ and $b$ in terms of $c$ using the Vandermonde transformation. We then choose the values of $c$ as the zeros of the shifted Legendre polynomial on the interval $[0,1]$.

\noindent Consider the Vandermonde matrix $V$ given as,
\renewcommand{\arraystretch}{1}
$$V=c^{j-1}_{i}=\begin{bmatrix}
  1 &c_{1} & c^{2}_{1} &\hdots&c^{s-1}_{1}\\
  1 & c_{2} & c^{2}_{2}&\hdots&c^{s-1}_{2} \\
  \vdots&\vdots&\vdots&\ddots&\vdots\\
   1&c_{s}&c^{2}_{s}&\hdots&c^{s-1}_{s}
\end{bmatrix}.$$
Multiply the symplectic condition \eqref{V1} of  Runge-Kutta methods with matrix $V$ as follows,
\begin{equation}\label{O}
c^{k-1}_{i}(b_{i}a_{ij}+b_{j}a_{ji}-b_{i}b_{j})c^{l-1}_{j}=0,
\hspace{.2cm} \forall \hspace{.1cm} i,j,k,l=1,2,\hdots,s.
\end{equation}
For methods with two stages ($s=2$), we take $l,k=1,2$, and then take summation over $i$ and $j$ from $1$ to $2$.\\
For $l=1$, $k=1,$
\begin{equation}\label{P}
\sum_{i,j} b_{i}a_{ij}+\sum_{i,j} b_{j}a_{ji}-\sum_{i,j}
b_{i}b_{j}=0.
\end{equation}
For $l=1$, $k=2,$
\begin{equation}\label{Q}
\sum_{i,j}b_{i}c_{i}a_{ij}+\sum_{i,j}b_{j}a_{ji}c_{i}-\sum_{i,j}b_{i}c_{i}b_{j}=0.
\end{equation}
For $l=2$, $k=1,$
\begin{equation}\label{R}
\sum_{i,j}b_{i}a_{ij}c_{j}+\sum_{i,j}b_{j}c_{j}a_{ji}-\sum_{i,j}b_{i}b_{j}c_{j}=0.
\end{equation}
For $l=2$, $k=2,$
\begin{equation}\label{S}
\sum_{i,j}b_{i}c_{i}a_{ij}c_{j}+\sum_{i,j}b_{j}c_{j}a_{ji}c_{i}-\sum_{i,j}b_{i}c_{i}b_{j}c_{j}=0.
\end{equation}
The following order two conditions must be satisfied.
\begin{equation}\label{T} \sum^{s}_{i=1}b_{i}=1, \hspace{.6cm}
\sum^{s}_{i=1}b_{i}c_{i}=\frac{1}{2}.
\end{equation}
Using equations \eqref{T} in equations \eqref{P}-\eqref{S} we have,
\begin{align*}&\displaystyle\sum_{i}b_{i}c_{i}&= \tfrac{1}{2},\\
&\displaystyle\sum_{i,j}b_{i}a_{ij}c_{j}+\displaystyle\sum_{i,j}b_{j}c_{j}a_{ji}&= \tfrac{1}{2}, \\
&\displaystyle\sum_{i,j}b_{i}c_{i}a_{ij}+\displaystyle\sum_{i,j}b_{j}a_{ji}c_{i}&= \tfrac{1}{2}, \\
&\displaystyle\sum_{i,j}b_{i}c_{i}a_{ij}c_{j} &= \tfrac{1}{8}.
\end{align*}
Consider the relation,
\[b_{i}(c_{i}-c_{1})= b_{i}c_{i}-b_{i}c_{1},\]
Take summation over $i$ from 1 to $s$ we get,
\[
\displaystyle\sum_{i}b_{i}(c_{i}-c_{1})=\displaystyle\sum_{i}b_{i}c_{i}-\displaystyle\sum_{i}b_{i}c_{1},\]
\[
b_{2}(c_{2}-c_{1})=\frac{1}{2}-c_{1},\]
\[
b_{2}=\frac{\tfrac{1}{2}-c_{1}}{c_{2}-c_{1}}.
\]
Similarly we can get,
\[b_{1} = \frac{\frac{1}{2}-c_{2}}{c_{1}-c_{2}}.\]
Now consider the relation,
\[b_{i}(c_{i}-c_{1})a_{ij}(c_{j}-c_{1})=b_{i}c_{i}a_{ij}c_{j}-b_{i}c_{i}a_{ij}c_{1}-b_{i}a_{ij}c_{j}c_{1}+b_{i}a_{ij}c_{1}c_{1}.   \]
Take summation over $i$ and $j$, and use previous equations we get,
\[ a_{22}=\frac{\frac{1}{8}-\frac{c_{1}}{3}-\frac{c_{1}}{6}+\frac{c_{1}c_{1}}{2}}{b_{2}(c_{2}-c_{1})(c_{2}-c_{1})}. \]
Similarly we get,
\[ a_{11}=\frac{\frac{1}{8}-\frac{c_{2}}{3}-\frac{c_{2}}{6}+\frac{c_{2}c_{2}}{2}}{b_{1}(c_{1}-c_{2})(c_{1}-c_{2})}, \]
\[ a_{21}=\frac{\frac{1}{8}-\frac{c_{2}}{3}-\frac{c_{1}}{6}+\frac{c_{1}c_{2}}{2}}{b_{2}(c_{2}-c_{1})(c_{1}-c_{2})}, \]
\[ a_{12}=\frac{\frac{1}{8}-\frac{c_{1}}{3}-\frac{c_{2}}{6}+\frac{c_{1}c_{2}}{2}}{b_{1}(c_{1}-c_{2})(c_{2}-c_{1})}. \]
A class of Runge-Kutta methods can be found by suitably choosing $c_{1}$ and $c_{2}$. We avail three possibilities in this regard and all are based on zeros of the shifted Legendre polynomials $P_{s}^{*}$ on the interval $[0,1]$ where,
\[
P_{s}^{*}(x)=\frac{s!}{2s}\displaystyle \sum_{k=0}^{s}(-1)^{s-k}\left(\begin{array}{c}
s\\
k
\end{array}\right)\left(\begin{array}{c}
s+k\\
k
\end{array}\right)x^{k}.
\]
For Gauss methods, the abscissa $c_i$ are the zeros of the shifted Legendre polynomials $P_{s}^{*}$ on the interval $[0,1]$ and has an order $2s$. For Radau methods, the first step is to choose the abscissa $c_{1}=0$ or $c_{s}=1$ or both of them. The rest of the abscissa are  chosen such that, for Radau  \Rmnum{1} methods, the abscissa are the zeros of the polynomial $P_{s}^{*}(x)+P_{s-1}^{*}(x)$ of order $2s-1$ or, for Radau  \Rmnum{2} methods, the abscissa are the zeros of the polynomial $P_{s}^{*}(x)-P_{s-1}^{*}(x)$ of order $2s-1$. For Lobatto  \Rmnum{3} methods, the abscissa are the zeros of the polynomial $P_{s}^{*}(x)-P_{s-2}^{*}(x)$ of order $2s-2$. Thus we have the following symplectic methods,

\textbf{Gauss, s=2:}\\
\renewcommand{\arraystretch}{1.8}\label{V11}
$$\begin{tabular}{l|*{6}{c}r}\label{V11}
$\frac{1}{2}-\frac{\sqrt{3}}{6}$  & $\frac{1}{4}$ & $\frac{1}{4}-\frac{\sqrt{3}}{6}$\\
$\frac{1}{2}+\frac{\sqrt{3}}{6}$  & $\frac{1}{4}+\frac{\sqrt{3}}{6}$ & $\frac{1}{4}$\\
\hline & $\frac{1}{2}$ & $\frac{1}{2}$ \label{V11}
\end{tabular}$$\\

\textbf{Radau  \Rmnum{1}, s=2:}\\
 \renewcommand{\arraystretch}{1.6}
$$\begin{tabular}{l|*{6}{c}r}
$0$ &$\frac{1}{8}$&$\frac{-1}{8}$\\
$\frac{2}{3}$ &$\frac{7}{24}$ & $\frac{3}{8}$\\
\hline & $\frac{1}{4}$&$\frac{3}{4}$
\end{tabular}$$\\

\textbf{Radau \Rmnum{2}, s=2:} \\
 \renewcommand{\arraystretch}{1.6}
$$\begin{tabular}{l|*{6}{c}r}
$\frac{1}{3}$ &$\frac{3}{8}$&$\frac{-1}{24}$\\
$1$ &$\frac{7}{8}$ & $\frac{1}{8}$\\
\hline & $\frac{3}{4}$&$\frac{1}{4}$
\end{tabular}$$
Similarly, we can construct methods with more stages and higher order.

\section{Construction of first integrals and their numerical preservation}\label{cfi}
\begin{figure}[H]
\centering
\minipage{0.5\textwidth}
\includegraphics[width=7cm,height=5.5cm]{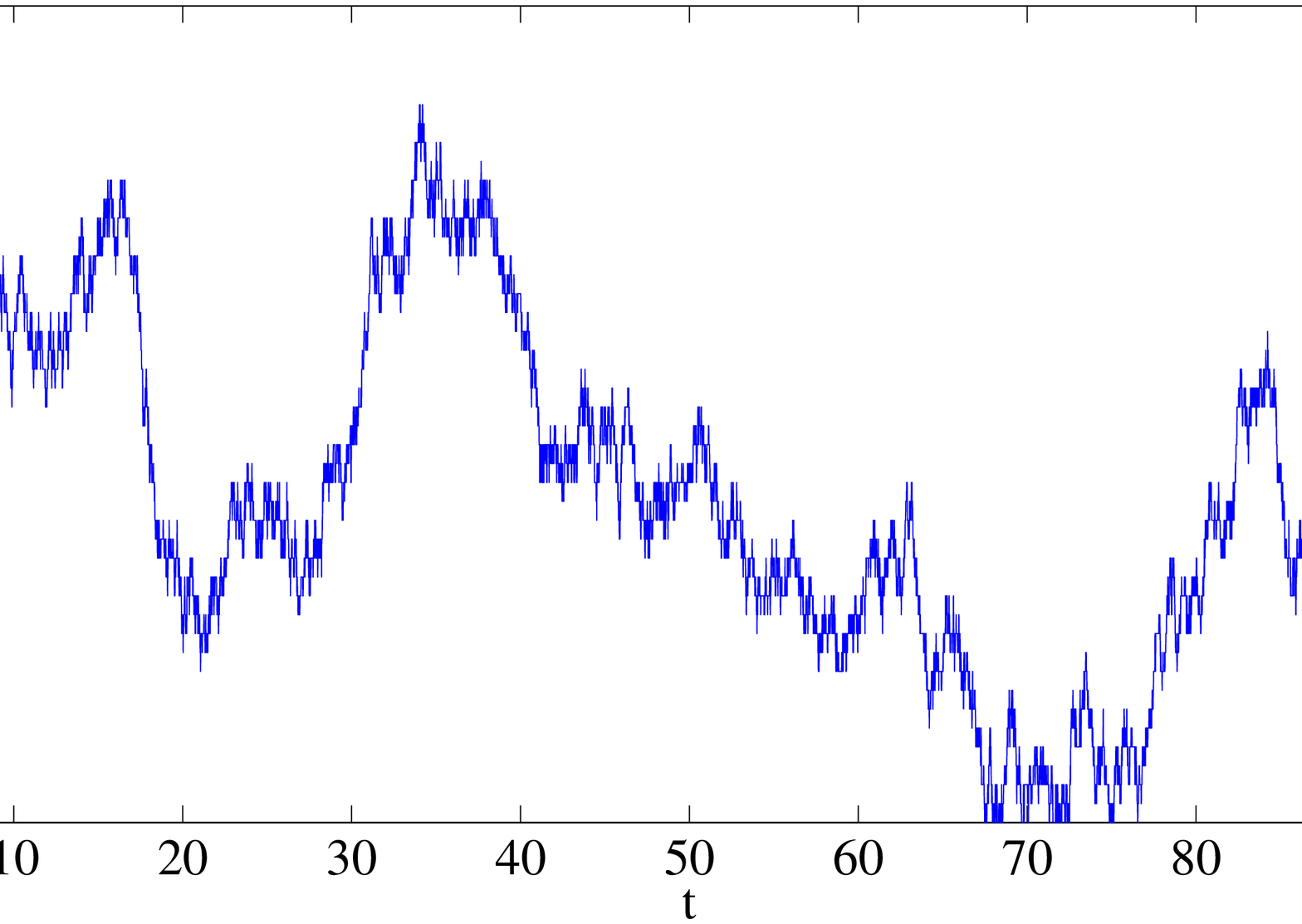}
\setlength\abovecaptionskip{-0.5cm}
\centering\caption{Error in integral $I_{2}$}\label{I2error}
\endminipage\hfill
\minipage{0.45\textwidth}
\includegraphics[width=7cm,height=5.5cm]{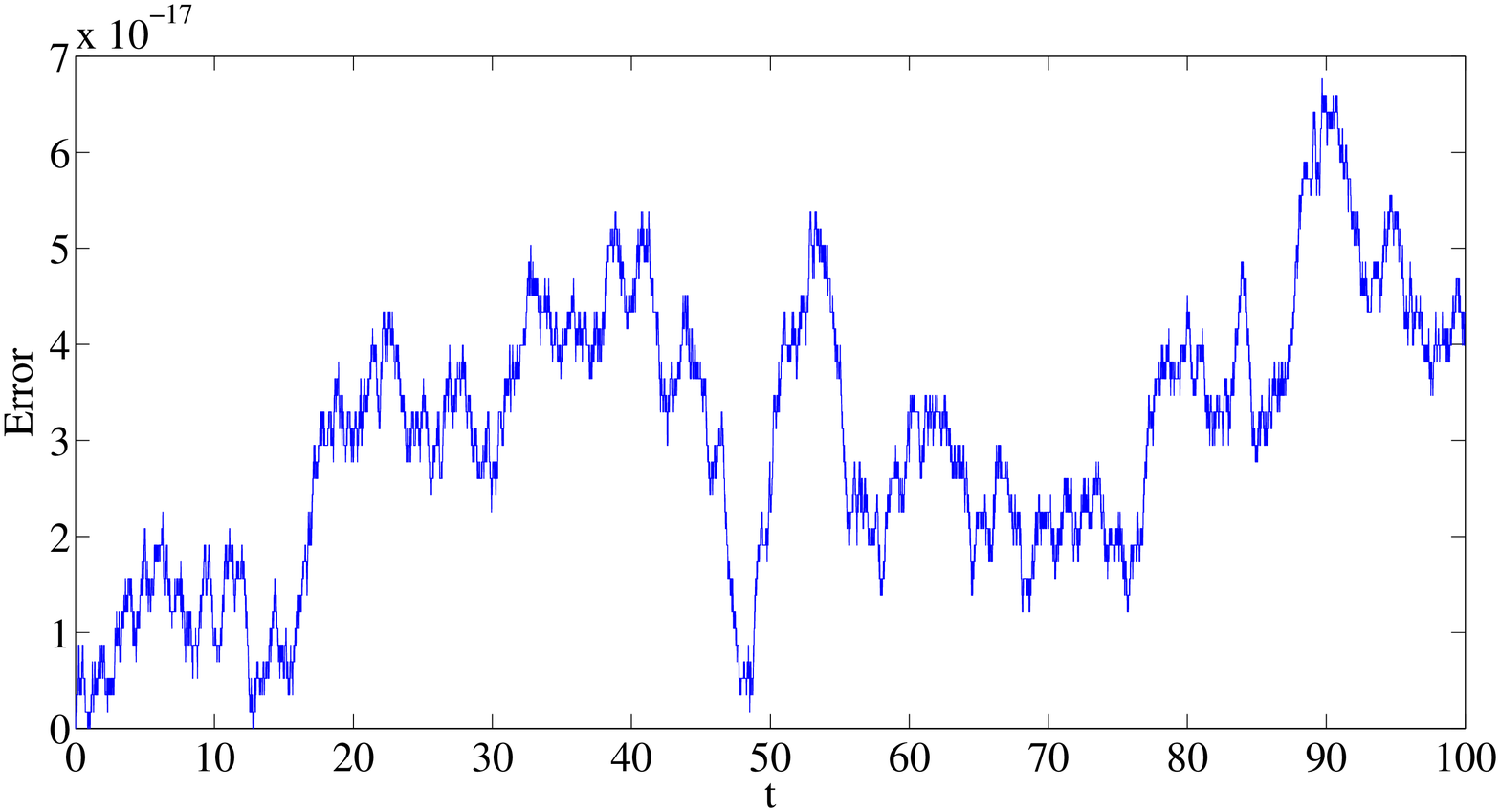}
\setlength\abovecaptionskip{-0.5cm}
\centering\caption{Error in integral $I_{3}$}\label{I3error}
\endminipage\hfill
\end{figure}

We construct first integrals of system of harmonic oscillators (both coupled and uncoupled) determined by the second order  ODE,
\begin{equation}\label{p1}
       y'' = -k^{2} y.
\end{equation}
We take different values of $k$ and $y$ as follows.\\
\noindent \textbf{Case I}:  {\textbf{(${{k^2=1}}$ and $y$ is real)}}\\\\
When $k^2=1$ and $y(t)$ is real valued, \eqref{p1} becomes one-dimensional harmonic oscillator equation
\begin{equation}\label{p2}
y''=-y,
\end{equation}
which possesses the standard Lagrangian
\begin{equation}\label{p3}
L=\frac{y'^{2}}{2}-\frac{y^{2}}{2}.
\end{equation}
Taking the Lagrangian and inserting in \eqref{4}, yields the
following determining system of equations
\begin{eqnarray}\label{p4}
\begin{aligned}
   - \eta{y} + {\eta_{t}}{y'} + {(\eta_{y}-\frac{1}{2}\xi_{t})}{y'^{2}} - \frac{1}{2}{\xi_{y}}{y'^{3}}-\frac{1}{2}{\xi_{t}}{y^{2}}
    -\frac{1}{2}{\xi_{y}}{y^{2}}{y'} - B_{t} - {y'}{B_{y}}=0.
\end{aligned}
\end{eqnarray}
Comparing different powers of $y'$ we have a system of four partial
differential equations whose solution gives rise to
\begin{equation}\label{p5}
\begin{aligned}
\xi{(t,y)} &= c_{1}+ c_{2}~{\sin{2t}}+c_{3}~{\cos{2t}},\\
\eta{(t,y)} &= {(c_{2}~{\cos{2t}}-c_{3}~{\sin{2t}})y}+c_{4}~{\sin{t}}+c_{5}~{\cos{t}},\\
B(t,y) &= -{(c_{2}~{\sin{2t}}+c_{3}~{\cos{2t}}){y^{2}}}+{(c_{4}~{\cos{t}}-c_{5}~{\sin{t}})}{y}.
\end{aligned}
\end{equation}
We thus obtain the following 5-Noether symmetry generators
\begin{eqnarray}\label{p6}
\begin{aligned}
X_{1} &= \frac{\partial}{\partial{t}},\\
X_{2} &= {\sin{2t}}~\frac{\partial}{\partial{t}}+{y}~{\cos{2t}}~\frac{\partial}{\partial{y}},\\
X_{3} &= {\cos{2t}}~\frac{\partial}{\partial{t}}-{y}~{\sin{2t}}~\frac{\partial}{\partial{y}},\\
X_{4} &= {\cos{t}}~\frac{\partial}{\partial{y}},\\
X_{5} &= {\sin{t}}~\frac{\partial}{\partial{y}}.\\
\end{aligned}
\end{eqnarray}
Using the symmetries \eqref{p6} and the Lagrangian \eqref{p3} in the
Noether's theorem \eqref{6}, we get following first integrals,
\begin{eqnarray}\label{p7}
\begin{aligned}
I_{1}&= \frac{y'^{2}}{2}+\frac{y^{2}}{2},\\
I_{2}&= {y'}~{\cos{t}}+{y}~{\sin{t}},\\
I_{3}&={y'}~{\sin{t}}-{y}~{\cos{t}},\\
I_{4}&= -\frac{1}{2}{y'^{2}}~{\cos{2t}}-{y}{y'}~{\sin{2t}} +\frac{1}{2}{y^{2}}~{\cos{2t}},\\
I_{5}&= -\frac{1}{2}{y'^{2}}~{\sin{2t}}+{y}{y'}~{\cos{2t}} +\frac{1}{2}{y^{2}}~{\sin{2t}}.
\end{aligned}
\end{eqnarray}
Amongst these five first integrals only two are independent
\cite{MZ}. We numerically integrate  \eqref{p2} using order four
Gauss $s=2$ symplectic Runge-Kutta method which we refer from now on
as Gauss2. We take step-size $h=0.01$ and number of steps
$n=10,000$. By employing a symplectic integrator, we expect the
first integrals of the system to be preserved by the
numerical scheme and this is what we have achieved. We look at the
deviation of numerically evaluated first integral $I(y_n)$ from the
actual value of first integral $I(y_0)$. We calculate error by
taking difference of the first integral evaluated at an initial
value $I(y_0)$ with the value of the first integral evaluated at all
subsequent numerically approximated values $I(y_n)$ given by the
formula: Error $ = |I(y_n)-I(y_{0})|$. Figure \ref{I2error} and
Figure \ref{I3error} represent the absolute error in the integral
$I_2$ and $I_3$ respectively. It is clear from the figures that the
error is very small and bounded, depicting qualitatively
correct numerical results. Similar error behavior is obtained for
other first integrals.

\begin{figure}
\centering
\minipage{0.5 \textwidth}
\includegraphics[width=7cm,height=5.5cm]{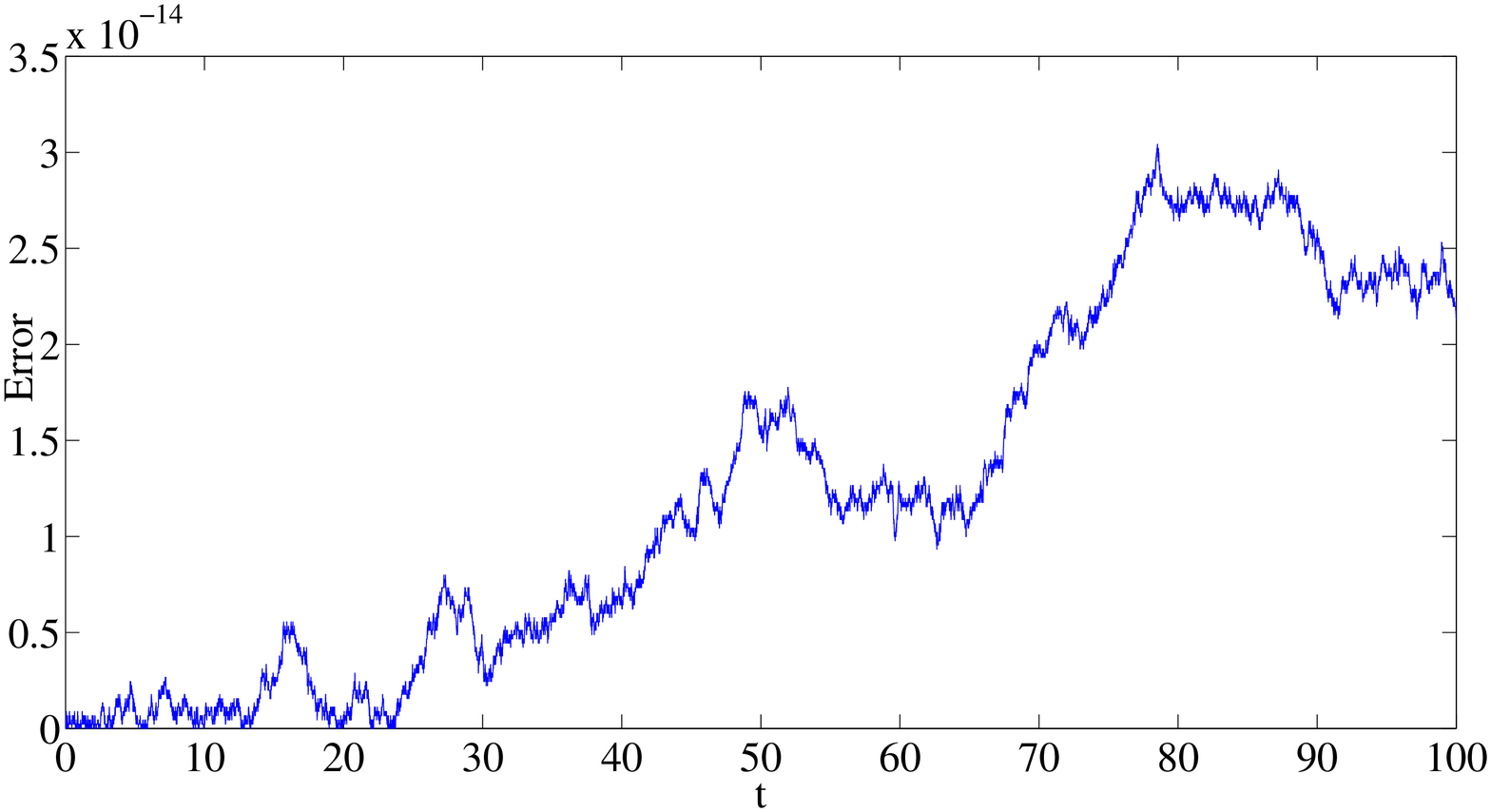}
\setlength\abovecaptionskip{-0.5cm}
\centering\caption{Error in integral $I_{2,1}$ }\label{pic3}
\endminipage\hfill
\minipage{0.45\textwidth}
\includegraphics[width=7cm,height=5.5cm]{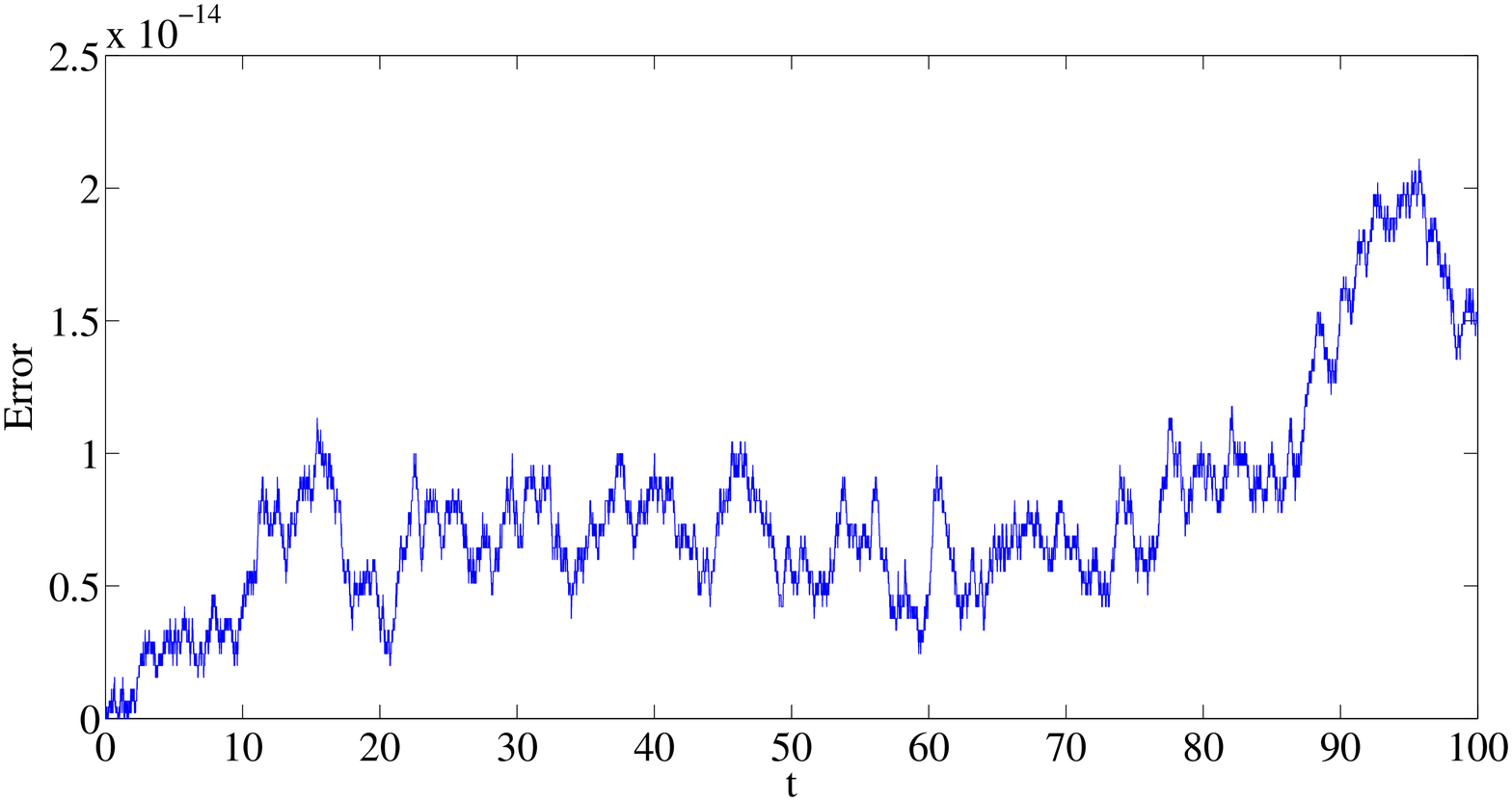}
\setlength\abovecaptionskip{-0.5cm}
\centering\caption{Error in integral $I_{2,2}$ }\label{pic4}
\endminipage\hfill
\vspace{1cm}
\minipage{0.5\textwidth}
\includegraphics[width=7cm,height=5.5cm]{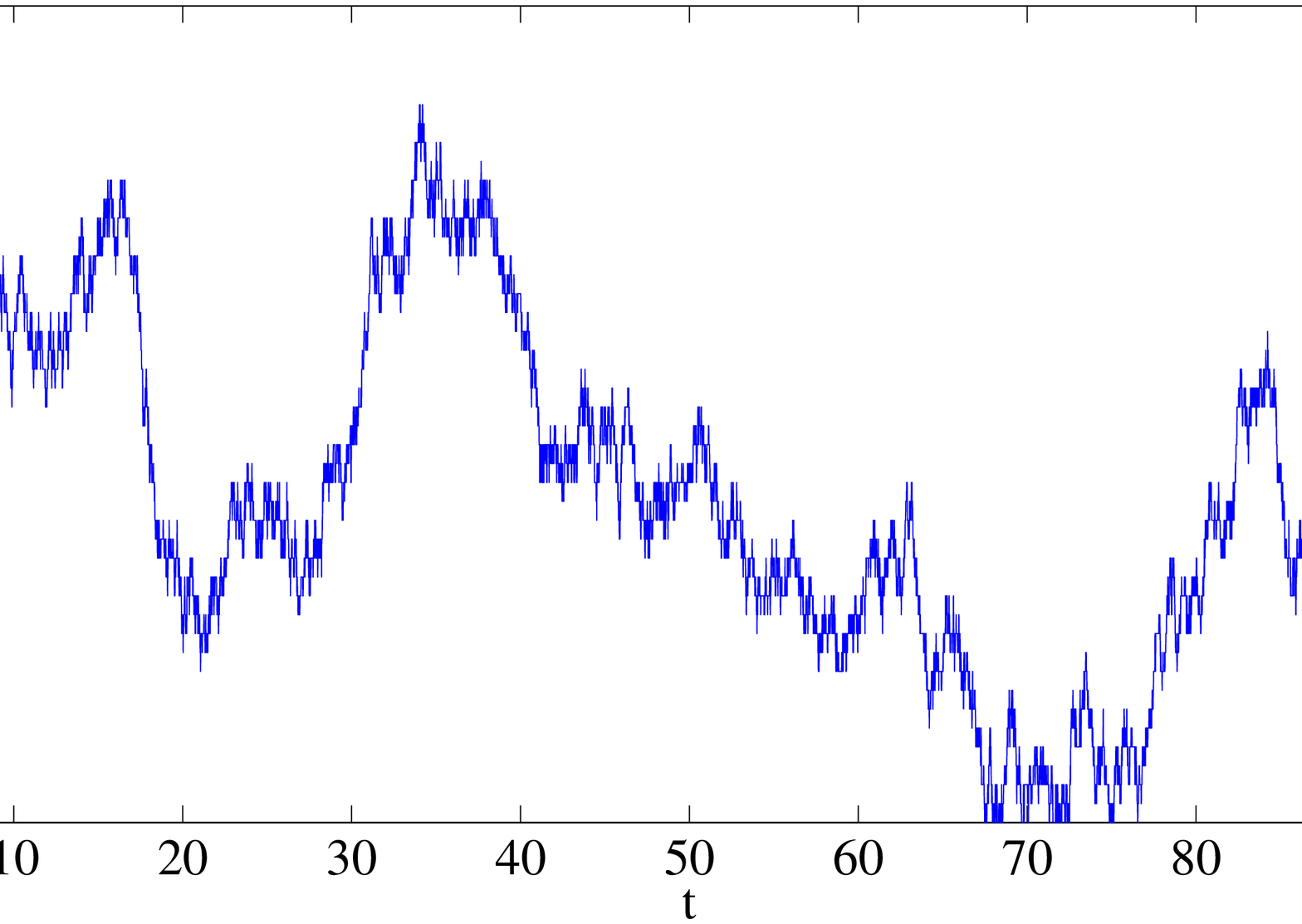}
\setlength\abovecaptionskip{-0.5cm}
\centering\caption{Error in integral $I_{4,1}$}\label{pic5}
\endminipage\hfill
\minipage{0.45\textwidth}
\includegraphics[width=7cm,height=5.5cm]{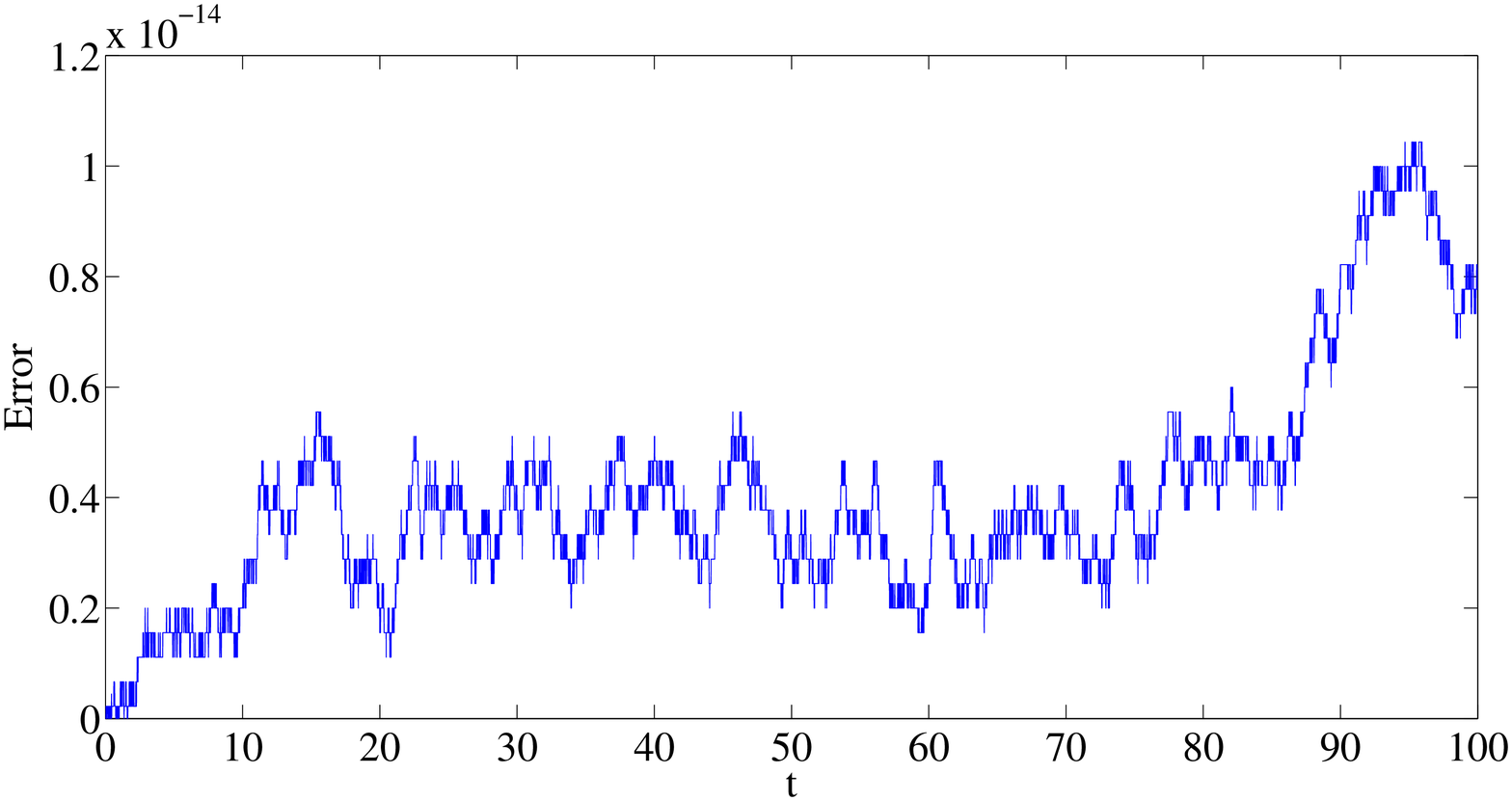}
\setlength\abovecaptionskip{-0.5cm}
\centering\caption{Error in integral $I_{4,2}$ }\label{pic6}
\endminipage\hfill
\end{figure}

%
%
%
\newpage
\noindent \textbf{Case II}: \textbf{($k^2=1$ and $y$ is complex)}\\ \\
When $k^2=1$ and $y(t)$ is a complex function $y=f+ig$ for $f$ and
$g$ being real functions of $t$, yields the following system of ODEs
\begin{equation}\label{p8}
\begin{aligned}
f''&=-f,\\
g''&=-g,
\end{aligned}
\end{equation}
which admits the Lagrangians
\begin{equation}\label{p9}
\begin{aligned}
L_{1}&=\frac{1}{2}(f'^{2}-g'^{2}-f^{2}+g^{2}),\\
L_{2}&=f'g'-fg.
\end{aligned}
\end{equation}
Using the Lagrangians \eqref{p9} in  \eqref{11} we get 9-Noether like operators
\begin{equation}\label{p10}
\begin{aligned}
X_{1}&=\frac{\partial}{\partial{t}},~~~~X_{2}=\sin{t}\frac{\partial}{\partial{f}},~~~~~X_{3}=\sin{t}\frac{\partial}{\partial{g}},
~~~~~X_{4}=\cos{t}\frac{\partial}{\partial{f}},~~~~X_{5}=\cos{t}\frac{\partial}{\partial{g}},\\
X_{6}&=\sin{2t}\frac{\partial}{\partial{t}}+f~{\cos{2t}}~\frac{\partial}{\partial{f}}+g~{\cos{2t}}~\frac{\partial}{\partial{g}},\\
X_{7}&=g~{\cos{2t}}~\frac{\partial}{\partial{f}}-f~{\cos{2t}}~\frac{\partial}{\partial{g}},\\
X_{8}&={\cos{2t}}~\frac{\partial}{\partial{t}}-f~{\sin{2t}}~\frac{\partial}{\partial{f}}-g~{\sin{2t}}~\frac{\partial}{\partial{g}},\\
X_{9}&=-g~{\sin{2t}}~\frac{\partial}{\partial{f}}+f~{\sin{2t}}~\frac{\partial}{\partial{g}}.
\end{aligned}
\end{equation}
\\
Invoking Eq. \eqref{12}, we obtain following invariants,
\begin{equation}\label{p11}
\begin{aligned}
I_{1,1}&=(f'^{2}-g'^{2}-f^{2}+g^{2})\sin{2t}-2(ff'-gg')\cos{2t},\\
I_{1,2}&=2(f'g'-fg)\sin2t-2(fg'+f'g)\cos{2t},\\
I_{2,1}&=(f'^{2}-g'^{2}-f^{2}+g^{2})\cos{2t}+2(f'f-g'g)\sin{2t},\\
I_{2,2}&=2(f'g'-fg)\cos{2t}+2(fg'+f'g)\sin{2t},\\
I_{3,1}&=-2f'\cos{t}-2f\sin{t},\\
I_{3,2}&=-2g'\cos{t}-2g\sin{t},\\
I_{4,1}&=-2f'\sin{t}+2f\cos{t},\\
I_{4,2}&=-2g'\sin{t}+2g\cos{t},\\
I_{5,1}&=f'^{2}-g'^{2}-f^{2}+g^{2},\\
I_{5,2}&=2f'g'+2fg.
\end{aligned}
\end{equation}
associated with Noether-like operators \eqref{p10}. System of
equations \eqref{p8} is integrated using Gauss2 method with stepsize
$h=0.01$ and number of steps $n=10,000$. The
absolute error in the first integrals $I_{2,1}$, $I_{2,2}$, $I_{4,1}$ and
$I_{4,2}$ is plotted in Figures \ref{pic3}, \ref{pic4}, \ref{pic5}
and \ref{pic6} respectively. Similar error behaviour is obtained for
$I_{1,1}$, $I_{1,2}$, $I_{3,1}$, $I_{3,2}$, $I_{5,1}$ and $I_{5,2}$.
We observe that the error does not grow out of bound which shows that
the numerical method is able to mimic the true qualitative feature
of the dynamical system.
\\

\begin{figure}[H]
\centering
\minipage{0.5\textwidth}
\includegraphics[width=7cm,height=5.5cm]{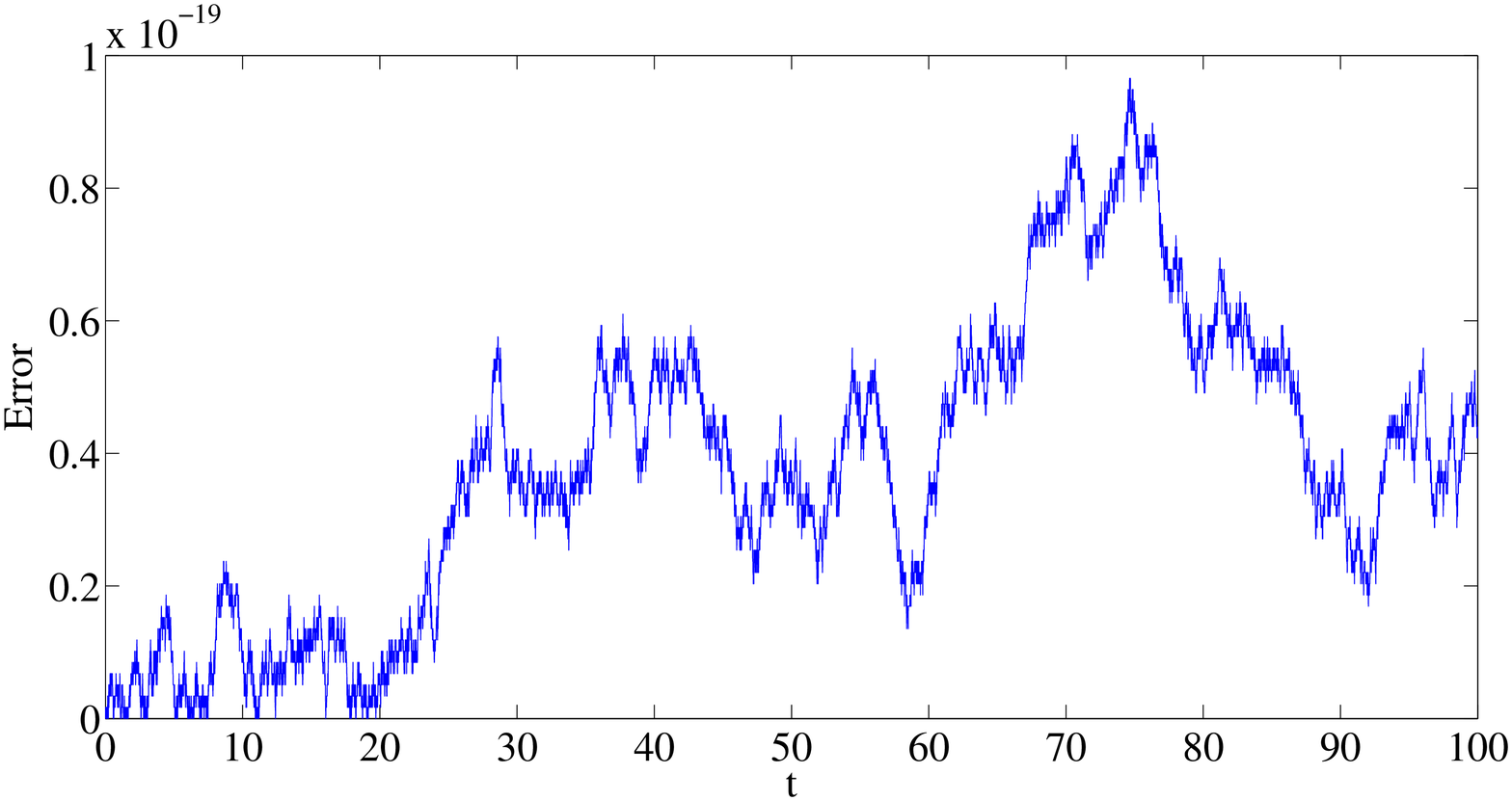}
\setlength\abovecaptionskip{-0.5cm}
\centering\caption{Error in integral $I_{1,1}$ }\label{pic7}
\endminipage\hfill
\minipage{0.45\textwidth}
\includegraphics[width=7cm,height=5.5cm]{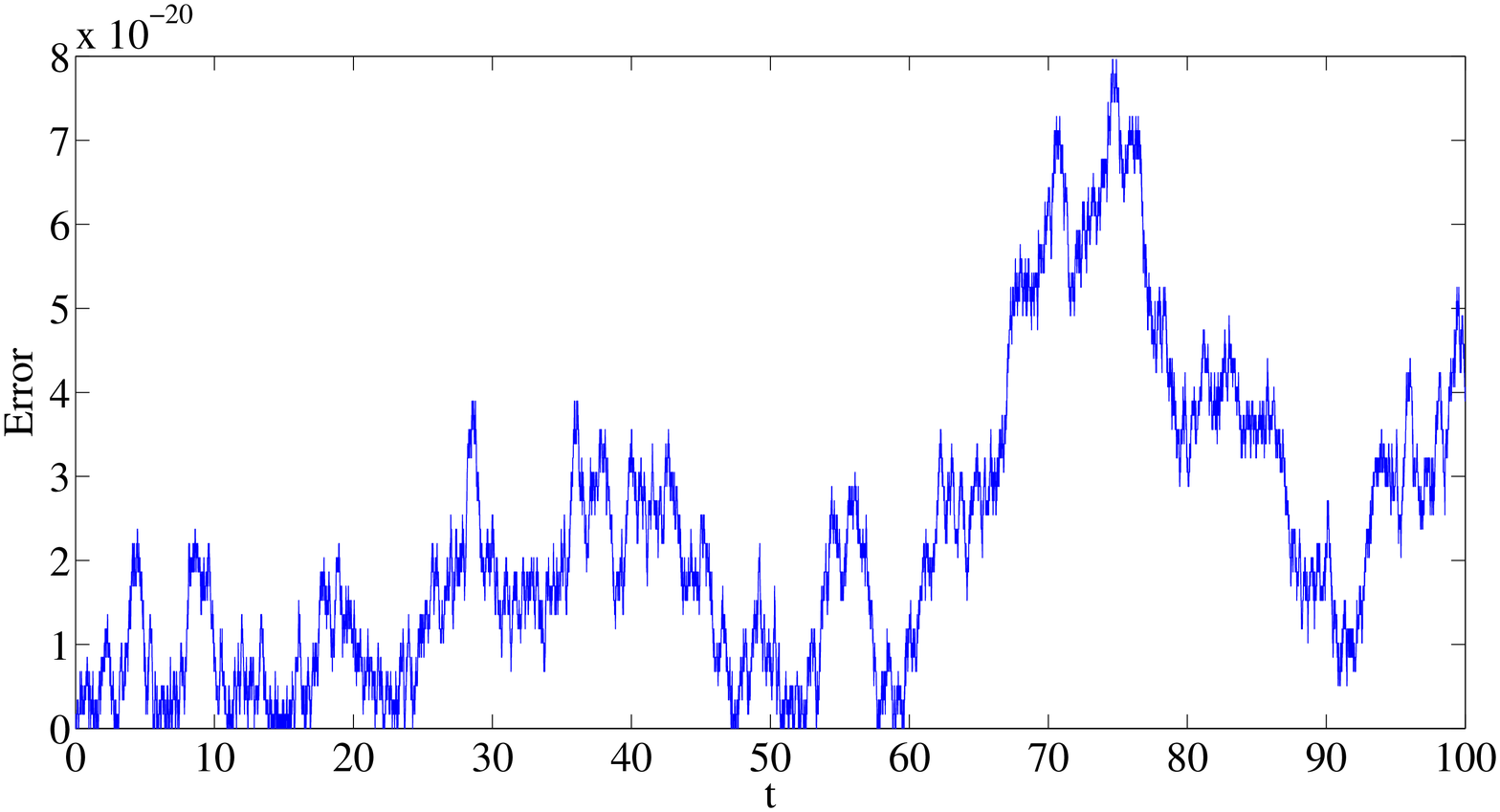}
\setlength\abovecaptionskip{-0.5cm}
\centering\caption{Error in integral $I_{1,2}$ }\label{pic8}
\endminipage\hfill
\vspace{2cm}
\minipage{0.5\textwidth}
\includegraphics[width=7cm,height=5.5cm]{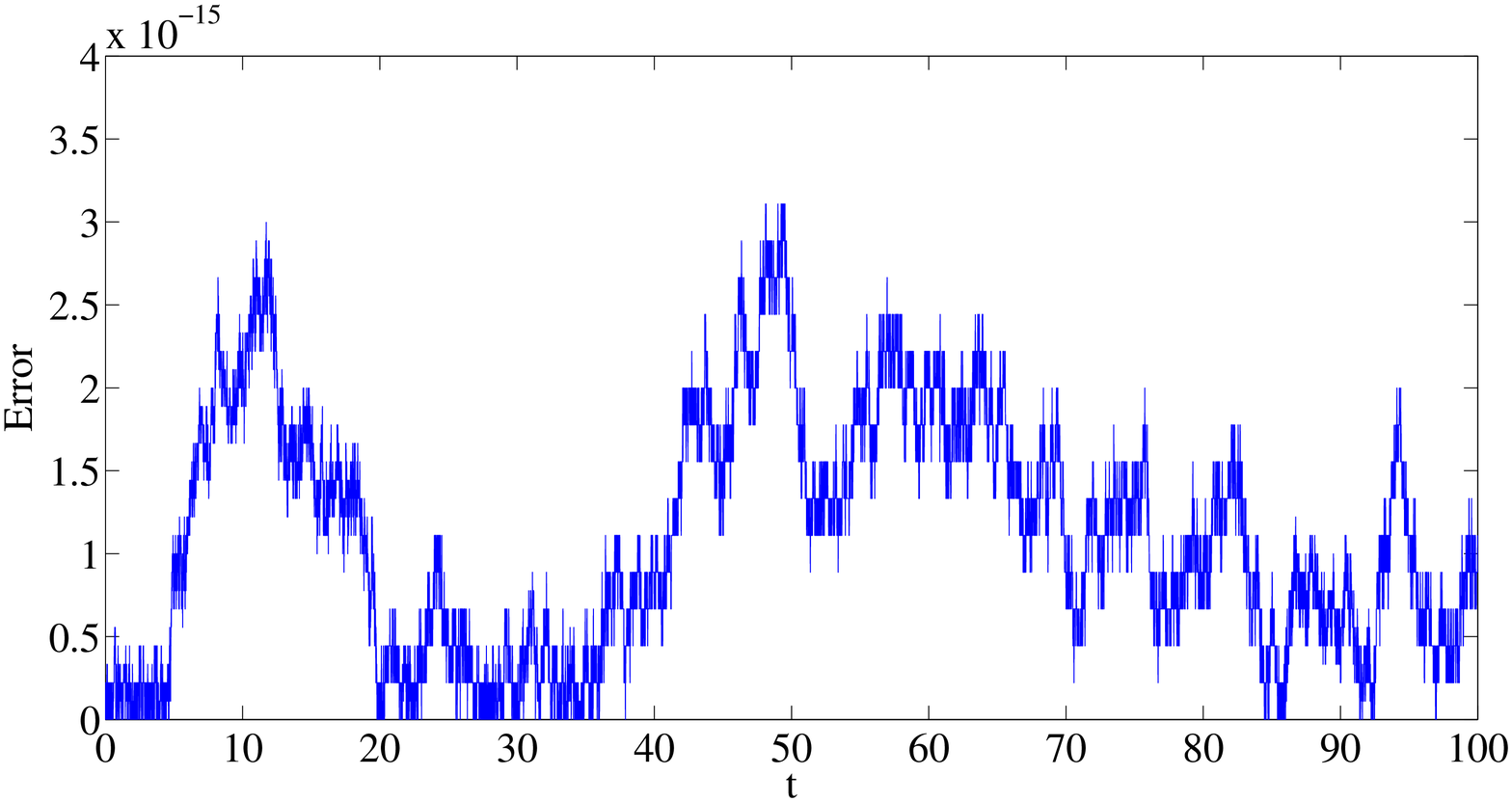}
\setlength\abovecaptionskip{-0.5cm}
\centering\caption{Error in integral $I_{3,1}$}\label{pic9}
\endminipage\hfill
\minipage{0.45\textwidth}
\includegraphics[width=7cm,height=5.5cm]{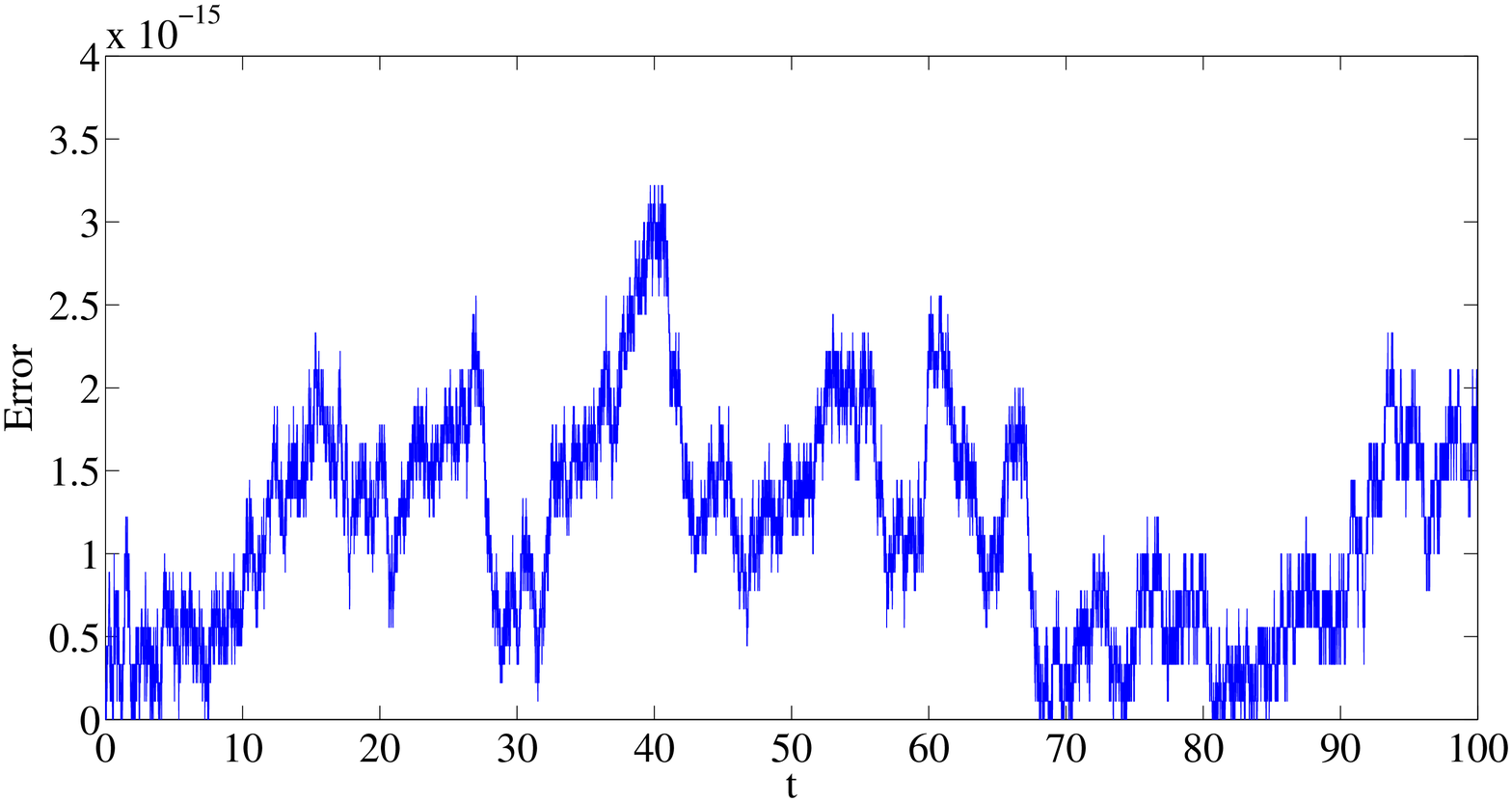}
\setlength\abovecaptionskip{-0.5cm}
\centering\caption{Error in integral $I_{3,2}$ }\label{pic10}
\endminipage\hfill
\end{figure}

\noindent\textbf{Case III}: \textbf{($k$ and $y$ are complex)}\\ \\
When $k$ and $y(t)$ are both complex, i.e., $k=\alpha_{1}+i\alpha_{2}$ and $y=f+ig$ for $f$, $g$, $\alpha_{1}$, and $\alpha_{2}$ being real,
 following coupled system of harmonic oscillators is obtained,\\
\begin{equation}\label{p12}
\begin{aligned}
f''&=-(\alpha_{1}^{2}-\alpha_{2}^{2})f+2\alpha_{1}\alpha_{2}g,\\
g''&=-(\alpha_{1}^{2}-\alpha_{2}^{2})g-2\alpha_{1}\alpha_{2}f,
\end{aligned}
\end{equation}
\\
which admits a pair of Lagrangians \cite{SA},
\begin{equation}\label{p13}
\begin{aligned}
L_{1}&=\frac{1}{2}(f'^{2}-g'^{2})- \frac{1}{2} (\alpha_{1}^2-\alpha_{2}^2)(f^{2}-g^{2})+2\alpha_{1}\alpha_{2}f g,\\
L_{2}&=f' g' - \alpha_{1}\alpha_{2} (f^{2}-g^{2}) -(\alpha_{1}^2-\alpha_{2}^2)f g.
\end{aligned}
\end{equation}
\\
The system \eqref{p12} admits following 9 Noether-like operators,
\begin{equation}\label{p14}
\begin{aligned}
X_{1}&=\frac{\partial}{\partial{t}},\\
X_{2}&=\sin({\alpha_{1}t})\cosh({\alpha_{2}t})\frac{\partial}{\partial{f}}+\cos({\alpha_{1}t})\sinh({\alpha_{2}t})\frac{\partial}{\partial{g}},\\
X_{3}&=\cos({\alpha_{1}t})\sinh({\alpha_{2}t})\frac{\partial}{\partial{f}}-\sin({\alpha_{1}t})\cosh({\alpha_{2}t})\frac{\partial}{\partial{g}},\\
X_{4}&=\cos({\alpha_{1}t})\cosh({\alpha_{2}t})\frac{\partial}{\partial{f}}-\sin({\alpha_{1}t})\sinh({\alpha_{2}t})\frac{\partial}{\partial{g}},\\
X_{5}&=-\sin({\alpha_{1}t})\sinh({\alpha_{2}t})\frac{\partial}{\partial{f}}-\cos({\alpha_{1}t})\cosh({\alpha_{2}t})\frac{\partial}{\partial{g}},\\
X_{6}&=\sin({2\alpha_{1}t})\cosh({2\alpha_{2}t})\frac{\partial}{\partial{t}}+\{(\alpha_{1}f-\alpha_{2}g)\cos({2\alpha_{1}t})\cosh({2\alpha_{2}t})+(\alpha_{1}g+\alpha_{2}f)\sin({2\alpha_{1}t})\sinh({2\alpha_{2}t})\}\frac{\partial}{\partial{f}}\\
&+\{(\alpha_{1}g+\alpha_{2}f)\cos({2\alpha_{1}t})\cosh({2\alpha_{2}t})-(\alpha_{1}f-\alpha_{2}g)\sin({2\alpha_{1}t})\sinh({2\alpha_{2}t})\}\frac{\partial}{\partial{g}},\\ X_{7}&=\cos({2\alpha_{1}t})\sinh({2\alpha_{2}t})\frac{\partial}{\partial{t}}+\{(\alpha_{1}g+\alpha_{2}f)\cos({2\alpha_{1}t})\cosh({2\alpha_{2}t})-(\alpha_{1}f-\alpha_{2}g)\sin({2\alpha_{1}t})\sinh({2\alpha_{2}t})\}\frac{\partial}{\partial{f}}\\
&-\{(\alpha_{1}f-\alpha_{2}g)\cos({2\alpha_{1}t})\cosh({2\alpha_{2}t})+(\alpha_{1}g+\alpha_{2}f)\sin({2\alpha_{1}t})\sinh({2\alpha_{2}t})\}\frac{\partial}{\partial{g}},\\ X_{8}&=\cos({2\alpha_{1}t})\cosh({2\alpha_{2}t})\frac{\partial}{\partial{t}}+\{(\alpha_{1}f-\alpha_{2}g)\sin({2\alpha_{1}t})\cosh({2\alpha_{2}t})-(\alpha_{1}g+\alpha_{2}f)\cos({2\alpha_{1}t})\sinh({2\alpha_{2}t})\}\frac{\partial}{\partial{f}}\\
&+\{(\alpha_{1}f-\alpha_{2}g)\cos({2\alpha_{1}t})\sinh({2\alpha_{2}t})+(\alpha_{1}g+\alpha_{2}f)\sin({2\alpha_{1}t})\cosh({2\alpha_{2}t})\}\frac{\partial}{\partial{g}},\\
X_{9}&=-\sin({2\alpha_{1}t})\sinh({2\alpha_{2}t})\frac{\partial}{\partial{t}}+\{(\alpha_{1}f-\alpha_{2}g)\cos({2\alpha_{1}t})\sinh({2\alpha_{2}t})+(\alpha_{1}g+\alpha_{2}f)\sin({2\alpha_{1}t})\cosh({2\alpha_{2}t})\}\frac{\partial}{\partial{f}}\\
&-\{(\alpha_{1}f-\alpha_{2}g)\sin({2\alpha_{1}t})\cosh({2\alpha_{2}t})-(\alpha_{1}g+\alpha_{2}f)\cos({2\alpha_{1}t})\sinh({2\alpha_{2}t})\}\frac{\partial}{\partial{g}}.\\\\
\end{aligned}
\end{equation}\\
Using Nother-like operators \eqref{p14} with pair of Lagrangians
\eqref{p13} in \eqref{12}, we obtain following ten first
integrals,
\begin{equation}\label{p15}
\begin{aligned}
I_{1,1}&=(\alpha_{1}^2-\alpha_{2}^{2})(f^{2}-g^{2})-4\alpha_{1}\alpha_{2}fg+f'^{2}-g'^{2},\\\\
I_{1,2}&=2(\alpha_{1}^2-\alpha_{2}^{2})fg+2\alpha_{1}\alpha_{2}(f^{2}-g^{2})+2f'g',\\\\
I_{2,1}&=f'\sin({\alpha_{1}t})\cosh({\alpha_{2}t})-g'\cos({\alpha_{1}t})\sinh({\alpha_{2}t})-(\alpha_{1}f-\alpha_{2}g)\cos({\alpha_{1}t})\cosh({\alpha_{2}t})\\
&-(\alpha_{1}g+\alpha_{2}f)\sin({\alpha_{1}t})\sinh({\alpha_{2}t}),\\\\
I_{2,2}&=g'\sin({\alpha_{1}t})\cosh({\alpha_{2}t})+f'\cos({\alpha_{1}t})\sinh({\alpha_{2}t})-(\alpha_{1}g+\alpha_{2}f)\cos({\alpha_{1}t})\cosh({\alpha_{2}t})\\
&+(\alpha_{1}f-\alpha_{2}g)\sin({\alpha_{1}t})\sinh({\alpha_{2}t}),\\\\
I_{3,1}&=f'\cos({\alpha_{1}t})\cosh({\alpha_{2}t})+g'\sin({\alpha_{1}t})\sinh({\alpha_{2}t})+(\alpha_{1}f-\alpha_{2}g)\sin({\alpha_{1}t})\cosh({\alpha_{2}t})\\
&-(\alpha_{1}g+\alpha_{2}f)\cos({\alpha_{1}t})\sinh({\alpha_{2}t}),\\\\
I_{3,2}&=g'\cos({\alpha_{1}t})\cosh({\alpha_{2}t})-f'\sin({\alpha_{1}t})\sinh({\alpha_{2}t})+(\alpha_{1}g+\alpha_{2}f)\sin({\alpha_{1}t})\cosh({\alpha_{2}t})\\
&+(\alpha_{1}f-\alpha_{2}g)\cos({\alpha_{1}t})\sinh({\alpha_{2}t}),\\\\
I_{4,1}&=\frac{1}{2}[\{(\alpha_{1}^{2}-\alpha_{2}^{2})(f^{2}-g^{2})-4\alpha_{1}\alpha_{2}fg-(f'^{2}-g'^{2})\}\sin(2{\alpha_{1}t})\cosh(2{\alpha_{2}t})\\
&-\{2\alpha_{1}\alpha_{2} (f^2-g^2)+2(\alpha_{1}^{2}-\alpha_{2}^{2})fg-2f'g'\}\cos({2 {\alpha_{1}}t})\sinh({2 \alpha_{2} t})]\\
&+\{\alpha_{1}(ff'-gg')-\alpha_{2}(fg'+gf')\}\cos({2\alpha_{1} t})\cosh({2\alpha_{2}t})\\
&+\{\alpha_{1}(fg'+gf')+\alpha_{2}(ff'-g'g)\}\sin({2\alpha_{1}t})\sinh({2\alpha_{2}t})\\\\
I_{4,2}&=\frac{1}{2}[\{(\alpha_{1}^{2}-\alpha_{2}^{2})(f^{2}-g^{2})-4\alpha_{1}\alpha_{2} fg-(f'^{2}-g'^{2})\}\cos({2\alpha_{1} t})\sinh({2\alpha_{2} t})\\
&+\{2\alpha_{1}\alpha_{2} (f^2-g^2)+2fg(\alpha_{1}^{2}-\alpha_{2}^{2})-2f'g'\}\sin({2\alpha_{1} t})\cosh({2\alpha_{2} t})]\\
&+[\{\alpha_{1}(fg'+f'g)+\alpha_{2}(ff'-gg')\}\cos({2\alpha_{1} t})\cosh({2\alpha_{2}t})]\\
&-[\{\alpha_{1}(ff'-gg')-\alpha_{2}(fg'+f'g)\}\sin({2\alpha_{1} t})\sinh({2\alpha_{2}t})]\\\\
I_{5,1}&=\frac{1}{2}[\{(\alpha_{1}^{2}-\alpha_{2}^{2})(f^{2}-g^{2})-4\alpha_{1}\alpha_{2} fg-(f'^{2}-g'^{2})\}\cos({2\alpha_{1} t})\cosh({2\alpha_{2}t})\\
&+\{2\alpha_{1}\alpha_{2}(f^2-g^2)+2(\alpha_{1}^{2}-\alpha_{2}^{2})fg-2f'g'\}\sin({2\alpha_{1} t})\sinh({2\alpha_{2} t})]\\
&+\{\alpha_{1}(fg'+gf')+\alpha_{2}(ff'-gg')\}\cos({2\alpha_{1} t})\sinh({2\alpha_{2} t})\\
&-\{\alpha_{1}(ff'-gg')-\alpha_{2}(fg'+gf')\}\sin({2\alpha_{1} t})\cosh({2\alpha_{2} t})\\\\
I_{5,2}&=\frac{1}{2}[\{-(\alpha_{1}^{2}-\alpha_{2}^{2})(f^{2}-g^{2})+4\alpha_{1}\alpha_{2}fg+(f'^{2}-g'^{2})\}\sin({2\alpha_{1}t})\sinh({2\alpha_{2} t})\\
&+\{2\alpha_{1}\alpha_{2} (f^2-g^2)+2fg(\alpha_{1}^{2}-\alpha_{2}^{2})-2f'g'\}\cos({2\alpha_{1} t})\cosh({2\alpha_{2} t})]\\
&-[\{\alpha_{1}(fg'+f'g)-\alpha_{2}(ff'-gg')\}\sin({2\alpha_{1}t})\cosh({2\alpha_{2} t})]\\
&-[\{\alpha_{1}(ff'-gg')-\alpha_{2}(fg'+f'g)\}\cos({2\alpha_{1}t})\sinh({2\alpha_{2}t})]\\\\
\end{aligned}
\end{equation}\\
Gauss2 method is again used to integrate \eqref{p12} with stepsize $h=0.01$ and number of steps $n=10,000$. The absolute error in the first integrals is calculated as before. The absolute error in integrals $I_{1,1}$, $I_{1,2}$, $I_{3,1}$ and $I_{3,2}$ is plotted in Figures \ref{pic7}, \ref{pic8}, \ref{pic9} and \ref{pic10} respectively, which remains bounded for long time. Similar error behaviour is obtained for $I_{2,1}$, $I_{2,2}$, $I_{4,1}$, $I_{4,2}$, $I_{5,1}$ and $I_{5,2}$. Symplectic Gauss2 method is able to preserve all first integrals obtained by employing complex symmetry analysis.
\newpage
\section{Conclusion}\label{con}
First integrals of the dynamical system $y''=-k^2 y$ are obtained
via classical Noether approach and complex symmetry method. The
later approach yields invariant energy as a particular example,
 that is stored in both oscillators. Since these first integrals are quadratic in nature,
 symplectic Runge--Kutta method, whose construction is also given in this paper,
 has successfully been applied to the system and numerical preservation of these
 first integrals have been obtained. Interestingly, the numerical method presented in this paper was able
  to preserve energy of the single oscillator as well as the energy stored in the
  pair of coupled oscillators that arise from complex Noether approach. The error in the
   first integrals remain bounded for long time which would not have been possible if we have employed non-symplectic integrators.


\begin{thebibliography}{99}


\bibitem{AQ} Ali S, Mahomed FM, Qadir A, Complex Lie symmetries for variational problems.  J. Nonlinear Math. Phys. 2008;5:25-35.

\bibitem{SAQ} Ali S, Mahomed FM, Qadir A, Complex Lie symmetries for scalar second-order ordinary differential equations. Nonlinear Anal.: Real World Appl. 2009;10:3335-3344.


\bibitem{BJ} Burrage K, Butcher JC, Stability criteria for implicit Runge-Kutta methods, SIAM J. Numer. Anal. 1979;16: 46-57.


\bibitem{rf3} Butcher JC, A history of Runge-Kutta methods, Appl. Numer. Math. 1996;20:247-260.



\bibitem{Stability_cooper} Cooper GJ, Stability of Runge-Kutta methods for trajectory problems,  IMA J. Numer. Anal. 1987;7:1-13.


\bibitem{UF} Farooq MU, Ali S, Mahomad FM, Two-dimensional systems that arise from the Noether classification of Lagrangians on the line, Appl. Math. Comput. 2011:6959-6973.

\bibitem{SA} Farooq MU, Ali S, Qadir A, Invariants of two-dimensional systems via complex Lagrangians with applications, Commun. Nonlinear Sci. Numer. Simul. 2011;16:1804-1810.

\bibitem{rf7} Habib Y, Long-term behaviour of G-symplectic methods, PhD Thesis, The University of Auckland, 2010.


\bibitem{rf4} Hairer E, Lubich C, Wanner G, Geometric Numerical Integration: Structure-Preserving Algorithms for Ordinary Differential Equations, second ed., Springer, 2005.


\bibitem{Kara} Ibragimov NH, Kara AH, Mahomed FM, Lie-Backlund and Noether
symmetries with applications, Nonlinear Dynamics 1998;15:115-136.


\bibitem{Ibragimov} Ibragimov NH, Elementary Lie Group Analysis and Ordinary Differential Equations, John Wiley and Sons, Chichester, UK, 1999.


\bibitem{SRK_LAS}  Lasagni FM, Canonical Runge-Kutta methods, ZAMP 1988;39:952-953.


\bibitem{Leach} Leach PGL, Applications of the Lie theory of extended groups in Hamiltonian mechanics: the oscillator and the Kepler problem, J. Aust. Math. Soc. 1981;23:173-186.


\bibitem{SL1} Lie S, Theorie der transformationsgruppen, Teubner, Leipzig, Germany, 1888.

\bibitem{SL2} Lie S, Vorlesungen \"{u}ber differentialgleichungen mit bekannten infinitesimalen transformationen, Teubner: Leipzig, Germany, 1891.

\bibitem{MZ} Lutzky M, Symmetry groups and conserved quantities for the harmonic oscillator, J. Phys. A: Math. Gen. 1978;11.

\bibitem{Naz} Naz R, Freire IL, Naeem I, Comparison of different approaches to construct first integrals for ordinary differential equations, Abstr. Appl. Anal., Hindawi Publishing Corporation 2014;1-15.

\bibitem{EN} Noether E, Invariante Variationsprobleme, Nachrichten der Akademie der Wissenschaften in G{\"o}ttingen, Mathematisch-Physikalische Klasse  1918;2:235-257. English translation in Transport Theory and Statistical Physics 1971;1(3):186-207.

\bibitem{Olver} Olver PJ, Applications of Lie Groups to Differential Equations, second ed., Springer, 1986.


\bibitem{rf20} Sanz-Serna JM, Runge-Kutta schemes for Hamiltonian systems, BIT 1988;28:877-883.
\bibitem{Serna} Sanz-Serna JM, Calvo MP, Symplectic numerical methods for Hamiltonian problems, J. Mod. Phys. C 1993;4:385-392.


\bibitem{Sanz} Sanz-Serna JM, Calvo MP, Numerical Hamiltonian Problems, Chapman and Hal, first edition, 1994.
\bibitem{Stephani}
Stephani H, Differential Equations: Their Solution Using Symmetries, Cambridge, UK, 1989.
\bibitem{rf30} Sun G, A simple way of constructing symplectic Runge-Kutta methods,
J. Comput. Math 2000;18:61-68.
%
%
%
%
%






%
%









\end{thebibliography}
\end{document}